\documentclass[12pt]{article}
\usepackage{amstext}
\usepackage{amsmath}
\usepackage{amssymb}
\numberwithin{equation}{section}
\begin{document}
\begin{center}
{\LARGE {\bf {\sf
Aspects of a new class of braid matrices: roots of unity and
hyperelliptic $q$ for triangularity, L-algebra,link-invariants,
noncommutative spaces  }}}
\\[0.8cm]

{\large  A.Chakrabarti\footnote{chakra@cpht.polytechnique.fr}},

\begin{center}
{\em
Centre de Physique Th\'eorique\footnote{Laboratoire Propre
du CNRS UPR A.0014}, Ecole Polytechnique, 91128 Palaiseau Cedex, France.\\
}
\end{center}

\end{center}

\smallskip

\smallskip

\smallskip

\smallskip

\smallskip

\smallskip

\begin{abstract}
Various properties of a class of braid matrices, presented before, are studied
considering $N^2 \times N^2 (N=3,4,...)$ vector representations for 
two subclasses. For $q=1$ the
matrices are nontrivial. Triangularity $(\hat R^2 =I)$ corresponds to 
polynomial equations for
$q$, the solutions ranging from roots of unity to hyperelliptic 
functions. The algebras of $L-$
operators are studied. As a crucial feature one obtains $2N$ central, 
group-like, homogenous
quadratic functions of $L_{ij}$ constrained to equality among 
themselves by the $RLL$
equations. They are studied in detail for $N =3$ and are proportional 
to  $I$ for the fundamental
$3\times3$ representation and hence for all iterated coproducts. The 
implications are analysed
through a detailed study of the $9\times 9$ representation for $N=3$. 
The Turaev construction for
link invariants is adapted to our class. A skein relation is 
obtained. Noncommutative spaces
associated to our  class of $\hat R$ are constructed. The transfer 
matrix map is implemented, with
the $N=3$ case as example, for  an iterated construction of 
noncommutative coordinates starting
from an $(N-1)$ dimensional commutative base space. Further 
possibilities, such as multistate
statistical models, are indicated.
\end{abstract}
\vfill

\newpage
\pagestyle {plain}
\section {Introduction :}
A new class of braid matrices was presented in previous papers. The 
most convenient formulation
can be found in Sec.$3$ of Ref.$1$.This is based on two previous 
works [$2,3$]. For ready
reference we summarize below the essential features. In succeeding 
sections different properties
of such braid matrices will be studied. Remarkable new aspects will 
be encountered. We will
$\it always$ be concerned with $N^2\times N^2$ vector representations 
of braid matrices
($N=3,4,...$).

For proper appreciation one should start by noting explicitly the 
links $\it and $ the crucial
differences with the standard $SO_q(N)$ and $Sp_q(N)$ braid matrices. 
Our approach is consistently
via spectral resolutions i.e. in terms of projectors.

The Baxterized  braid matrices ( depending on  a spectral parameter 
$\theta$ ) satisfy
\begin{equation}
{\hat R}_{12}(\theta) {\hat R}_{23}(\theta + \theta') {\hat R}_{12}(\theta') =
  {\hat R}_{23}(\theta'){\hat R}_{12}(\theta +\theta') {\hat R}_{23}(\theta)
\end{equation}
where
$${\hat R}_{12} =\hat R \otimes I_N, \qquad \hat R_{23}= I_N \otimes \hat R$$
and $I_N$ is the $N\times N$ identity matrix.

For the standard $SO_q(N)$ and $Sp_q(N)$ cases ( see sources cited in 
Sec.$2$ of Ref.$1$ ) one has
\begin{equation}
\hat R (\theta) = P_+ +v(\theta) P_- + w(\theta) P_0
\end{equation}
where the projectors satisfy
\begin{equation}
P_i P_j = \delta _{ij} P_i, \qquad P_{+} + P_{-} + P_{0} = I_{N^2}
\end{equation}
All $\theta$ - dependence is in $v(\theta)$ and $w(\theta)$. The 
projectors depend only on $q$.

The  elegant  canonical formulation ( Sec.$2$, Ref.$1$ ) gives the 
follwing $v(\theta)$ and
$w(\theta)$ where
$$q=exp h$$
One has for each case (independently of $N$)
\begin{equation}
v(\theta) = \frac { sinh(h-\theta)}{ sinh(h+\theta)}
\end{equation}
and two solutions for $w(\theta)$.

$SO_q(2n+1)$:
\begin{equation}
w(\theta)=\frac{cosh((n+\frac{1}{2})h -\theta)}{cosh((n+\frac{1}{2})h 
+\theta)};\qquad
\frac{sinh((n-\frac{1}{2})h -\theta)}{sinh((n-\frac{1}{2})h +\theta)}v(\theta)
\end{equation}

$SO_q(2n)$:
\begin{equation}
w(\theta)=\frac{cosh(nh -\theta)}{cosh(nh +\theta)};\qquad
\frac{sinh((n-1)h -\theta)}{sinh((n-1)h +\theta)}v(\theta)
\end{equation}

$Sp_q(2n)$:
\begin{equation}
w(\theta)=\frac{sinh((n+1)h -\theta)}{sinh((n+1)h +\theta)};\qquad
\frac{cosh(nh -\theta)}{cosh(nh +\theta)}v(\theta)
\end{equation}

In contrast, for our new class (Sec.$3$, Ref.$1$), ${\it conserving}$ 
the projectors but changing
the coefficients, we have
\begin{equation}
v(\theta)=1, \qquad w(\theta) = \frac{sinh(\eta -\theta)}{sinh(\eta +\theta)}
\end{equation}
where
\begin{equation}
e^{\eta} + e^{-\eta} = ([N-\epsilon] +\epsilon) = \frac 
{q^{N-\epsilon} -q^{-N +\epsilon}}{
q-q^{-1}} +\epsilon
\end{equation}
and $\epsilon =\pm1$ when the $P_i$ are those for $SO_q(N)$ and 
$Sp_q(N)$ respectively.
(An overall ambiguity of sign for the right side of $(1.9)$ has been 
fixed to assure real $\eta$
for real $q$. This will be maintained throughout, though complex $q$ 
will be considered later.)

We adopt  the following notations for our $\hat R (\theta)$ when 
$(1.8)$ is implemented:

$(a)$: When $P_i$ are those for $SO_q$ our $\hat{ R}(\theta)$ is of 
type $\hat{o}(N)$

$(b)$:  When $P_i$ are those for $Sp_q$ our $\hat{ R}(\theta)$ is of 
type $\hat{p}(N)$ .

This is to signal the provenance of the projectors and also at the 
same time the fact that the
coefficients $(1.8)$ often lead to startlingly different properties 
as compared to the standard
cases ( from $(1.4)$ to $(1.7)$ ).

For $(1.8)$ we obtain
\begin{equation}
\hat R (\theta) = P_+ + P_- + \frac{sinh(\eta -\theta)}{sinh(\eta 
+\theta)} P_0 \qquad = I +
\biggl(\frac{sinh(\eta -\theta)}{sinh(\eta +\theta)} -1\biggr) P_0
\end{equation}
For completeness we give $P_0$ explicitly $[4,5]$. Let the n-tuple 
$(\rho_1,\rho_2,...,\rho_N)$ be
defined as follows for the respective cases indicated:

\begin{equation}
  SO_q(2n+1): \qquad ( n-\frac{1}{2}, n -\frac{3}{2}, ... 
,\frac{1}{2},0,-\frac{1}{2},... ,
-n+\frac{1}{2} )
\end{equation}

\begin{equation}
  SO_q(2n): \qquad ( n-1, n -2, ... ,1,0,0,-1,... , -n+1 )
\end{equation}

\begin{equation}
  Sp_q(2n): \qquad ( n, n -1, ... ,1,-1,... , -n )
\end{equation}

Define correspondingly for
\begin{equation}
  SO_q(N): \qquad \epsilon =1, \qquad  (i= 1,...,N; \quad N= 2n, 2n+1 )
\end{equation}
\begin{equation}
  Sp_q(2n): \qquad \epsilon =1 \quad (i \leq n ), \quad \epsilon =-1 
\quad (i > n )
\qquad  (i= 1,...,2n )
\end{equation}
Set
  \begin{equation}
  i' = N+1 -i
\end{equation}
Then corresponding to $(1.11),(1.12)$, in notations most suitable for  us,
\begin{equation}
([N-1] +1) P_0 = \sum_{i,j =1}^{N}q^{(\rho_{i'} - \rho _{j})} 
(ij)\otimes (i'j') \quad \equiv P'_0
\end{equation}
and corresponding to $(1.13)$ and $(1.15)$
\begin{equation}
([N+1] -1) P_0 = \sum_{i,j =1}^{N}q^{(\rho_{i'} - \rho _{j})} 
\epsilon_{i} \epsilon_{j}(ij)\otimes
(i'j') \quad \equiv P'_0
\end{equation}
Here $(ij)$ denotes the $N\times N$ matrix with $1$ at (row-$i$, 
col.-$j$) and zero elsewhere.

These standard $P_0$ will be ${\it carried}$ ${\it over}$ to our 
$\hat o (2n+1),\hat o (2n),\hat  p
(2n)$. We do not indicate the $q$-dependence explicitly ( by denoting 
$\hat o _{q} (2n+1)$ for
example ) since for $q=1$ our constructions remain nontrivial. Our 
$\hat R$ matrices are $\it not$
$q$-deformations of a "classical" limit for a particular value of $q$ 
such as $1$. This is just one
of the remarkable features to be studied below.

For $N=3$, the $9\times 9$ projector $P_0$ is given by
\begin {equation}
(q+1+q^{-1})P_0 \equiv P'_0=\begin{vmatrix} 0 & 0 & 0 & 0 & 0 & 0 & 0 
& 0 & 0 \\ 0 & 0 & 0 & 0 & 0
& 0 & 0 & 0 & 0 \\ 0 & 0 & q^{-1} & 0 & q^{-\frac{1}{2}} & 0 & 1 & 0 
& 0 \\  0 & 0 & 0 & 0 & 0 & 0
& 0 & 0 & 0\\ 0 & 0 & q^{-\frac{1}{2}} & 0 & 1 & 0 & q^{\frac{1}{2}} 
& 0 & 0\\ 0 & 0 & 0 & 0 & 0 &
0 & 0 & 0 & 0\\ 0 & 0 & 1 & 0 & q^{\frac{1}{2}} & 0 & q & 0 & 0 \\ 0 
& 0 & 0 & 0 & 0 & 0 & 0 & 0 &
0 \\ 0 & 0 & 0 & 0 & 0 & 0 & 0 & 0 & 0 \end {vmatrix}
\end{equation}

This is the case that will be studied here extensively as the simplest example. For 
all $N$ a
basic feature is the proportionality of the rows with nonzero 
elements. This has important
consequences.
It remains here to display briefly the "pre-Baxterized" situation. 
Our canonical forms ensure,
among various aspects studied in Ref.$1$,
\begin{equation}
\hat R (-\theta) = (\hat R (\theta))^{-1},  \qquad    \hat R (0) = I_{N^2}
\end{equation}
The limits
$$ \theta \rightarrow \pm \infty, \qquad \hat R (\theta) \rightarrow 
{\hat R}^{\pm1} $$
satisfy the ( non-Baxterized, $\theta$-independent ) braid equation
\begin{equation}
{\hat R}_{12} {\hat R}_{23} {\hat R}_{12} =   {\hat R}_{23}{\hat 
R}_{12} {\hat R}_{23}
\end{equation}
where one can substitute  $ {\hat R}^{-1}$ for $ \hat R$.

For the standard cases
\begin{equation}
(SO_q (N)):{\hat R}^{\pm1} = P_+ - q^{\mp 2}P_- +q^{\mp N} P_0
\end{equation}
\begin{equation}
(Sp_q (N)):{\hat R}^{\pm1} = P_+ - q^{\mp 2}P_- - q^{\mp (N+2)} P_0
\end{equation}
satisfy $\it { cubic}$ equations. For our cases, with $\eta$ given by 
$(1.9)$, one has for all $N$
the $\it {quadratic}$ equation
\begin{equation}
(e^{\eta} \hat R) - (e^{\eta} \hat R)^{-1} = (e^{\eta} - e^{-\eta})I
\end{equation}
or
\begin{equation}
(\hat R -I) ( \hat R + e^{-2\eta}I) =0
\end{equation}

 From $(1.10)$, as $\theta \rightarrow \pm \infty $,
\begin{equation}
\hat R^{\pm1} = I - ( 1+ e^{\mp 2 \eta}) P_0 = I -  e^{\mp  \eta}( 
e^{ \eta}+ e^{-\eta}) P_0
\qquad = I -   e^{\mp \eta} P'_0
\end{equation}

The last equation follows from $(1.9),(1.17)$ and $(1.18)$. Note that 
though $P'_0$ is {\it {not}}
a projector in $(1.26)$, $\hat R$ is inverted by inverting the 
coefficient of $P'_0$ due to the
relation
$${P'_0}^2 = ( e^{\eta} +e^{-\eta}) P'_0 $$
Other properties of $\hat R$ will be introduced later as they become 
directly relevant.

\section { What $q$ for triangularity ? :}
A braid matrix for vector representation is called "triangular" if
\begin{equation}
            {\hat R}^2 =I
\end{equation}
For the standard  cases $(A,B,C,D)_q$ this is obtained trivially for 
$q=1$. This is well-known.
But for comparison with our case let  us  briefly indicate how this 
happens for $SO_q(N)$ and
$Sp_q(N)$.

For the projectors in $(1.2)$ and $(1.3)$ denote
\begin{equation}
(P_i)_{q=1} = \textsf{P}_i \qquad (i= +,-,0)
\end{equation}
and let
\begin{equation}
P = \sum_{i,j} (ij)\otimes (ji), \qquad P^2 =I
\end{equation}
(Acting on the left $P$ permutes specific rows and acting on the 
right the corresponding columns.
This evident feature is mentioned since it plays a crucial role below 
$(2.40)$.)

 From $(1.22),(1.23)$ substituing the known explicit forms $[4,5]$ of 
the projectors for $q=1$,
with upper and lower signs for the two cases respectively,
\begin{equation}
\hat R =\textsf{P}_+ - \textsf{P}_- \pm \textsf{P}_0 = \pm (I 
-2\textsf{P}_{\mp}) = P
\end{equation}
Hence
\begin{equation}
            {\hat R}^2 = P^2 =I, \qquad R = P\hat R = I
\end{equation}
for both cases.

For $GL_q(N)$ one  obtains the same result even more simply.

For  our  class $q=1$ gives a quite nontrivial situation, as 
emphasized already in Ref.$3$.
Denoting for all $N = (3,4,...)$
$$  {(\eta)}_{q=1} = \hat {\eta} $$
from $(1.9)$
\begin{equation}
e^{\hat {\eta}} +  e^{-\hat {\eta}} = N, \qquad  (1+  e^{\mp 2\hat 
{\eta}}) = \frac {2N}{N\pm
(N^2 -4)^\frac {1}{2}} \neq 2
\end{equation}
and from $(1.26)$
\begin{equation}
({\hat R}^{\pm1})_{q=1} = I - (1+  e^{\mp 2\hat {\eta}} ) \textsf 
{P}_0 , \qquad  {\hat R}^2 \neq 1
\end{equation}

The generalized Hecke  condition is now
\begin{equation}
(e^{\hat {\eta}} \hat R) - (e^{\hat {\eta}} \hat R)^{-1} = (e^{\hat 
{\eta}} - e^{-\hat {\eta}})I
\end{equation}
This cannot be conjugated to
\begin{equation}
(\hat R -I)(\hat R +I) =0
\end{equation}

For $(2.1)$ we need for our case
\begin{equation}
\eta =0
\end{equation}
when
$$ \hat R = I - 2 P_0, \qquad \hat R^2 = I +4({P_0}^2 -P_0)= I$$
Hence from $(1.9)$ for $\hat {o}(N)$ and  $\hat {p}(N)$ respectively
\begin{equation}
[N\mp 1]\pm 1 = e^{\eta} + e^{-\eta} = 2
\end{equation}
or respectively,
\begin{equation}
(A): \quad  q^{N-2}+ q^{N-4} + ...+ q^{-N+4} + q^{-N+2} = 1, \qquad (N=3,4,...)
\end{equation}
\begin{equation}
(B): \quad  q^{N}+ q^{N-2} + ...+ q^{-N+2} + q^{-N} = 3, \qquad (N=4,6,...)
\end{equation}
  The degrees of the polynomials can be lowered by changing variables 
as follows for the different
cases. To start with $(A)$ is divided into two subclasses.
$$(A)_1: N=2n+2 \qquad ( n=1,2,...)$$
Set
\begin{equation}
p=q^2,\qquad Y= p + p^{-1}
\end{equation}
when
\begin{equation}
q=\pm p^{\frac {1}{2}} = \pm \frac {1}{\sqrt 2} \biggl ( Y + \sqrt 
{Y^2 -4} \biggr)^{\frac {1}{2}}
\end{equation}
 From $(2.12)$,
\begin{equation}
S_n \equiv (p^n +p^{-n}) +(p^{n-1} +p^{-n+1}) + ...+ (p^2 +p^{-2}) + 
(p +p^{-1})  =0
\end{equation}
Note that for $(A)_1$ the right side cancels with $p^0 =1$. Now to 
express $S_n$ in terms of $Y$
implementing
\begin{equation}
S_{n+1} = YS_n - S_{n-1} +Y -2
\end{equation}
one obtains finally
$$S_n =Y^n + Y^{n-1} - \binom{n-1}{1} Y^{n-2} - \binom{n-2}{1} 
Y^{n-3} + \binom{n-2}{2} Y^{n-4}
+\binom{n-3}{2} Y^{n-5}+ ...$$
\begin{equation}
...+(-1)^r \biggl(\binom{n-r}{r} Y^{n-2r} +\binom{n-r-1}{r} 
Y^{n-2r-1}\biggr)+...+c_1 Y +c_0
\end{equation}
where
$$ c_1 = (-1)^{s-1}s, \qquad (n=2s, r=s-1); \qquad c_1 = 
(-1)^{s}(s+1), \qquad (n=2s+1, r=s)$$
\begin{equation}
c_0 = -2  \qquad (n= 2+4m,3+4m;\quad m=0,1,2,...); \qquad c_0 =0 
\qquad (n\neq  2+4m,3+4m)
\end{equation}

Explicitly
\begin{eqnarray}
\nonumber
S_1&=&Y\\
\nonumber
  S_2 &=& Y^2 + Y -2 = (Y-1)(Y+2)\\
\nonumber
S_3 &=& Y^3 +Y^2 -2Y -2 = (Y+1)(Y^2 -2)\\
\nonumber
S_4&=& Y^4+Y^3 -3Y^2 -2Y\\
\nonumber
  S_5& = &Y^5 +Y^4-4Y^3- 3Y^2 +3Y\\
\nonumber
S_6&= &Y^6+Y^5-5Y^4 -4Y^3 +6Y^2 +3Y -2\\
\nonumber
S_7&=& Y^7+Y^6-6Y^5 -5Y^4 +10Y^3 +6Y^2 -4Y -2\\
\nonumber
S_8&=& Y^8+Y^7-7Y^6 -6Y^5 +15Y^4 +10Y^3 -10Y^2 -4Y \\
S_9&=& Y^9+Y^8-8Y^7 -7Y^6 +21Y^5+15Y^4 -20Y^3 -10Y^2 +5Y
\end{eqnarray}
and so on.
$$(A)_2: N=2m+1 \qquad ( m=1,2,...)$$
Set
\begin{equation}
z=q+q^{-1}, \qquad q^{\pm 1} =\frac{1}{2}(z \pm \sqrt {z^2 -4})
\end{equation}
Retaining only the odd powers in $(2.18),(2.20)$ and adapting 
notations one obtains
\begin{eqnarray}
\nonumber
  \Sigma _{(2m-1)}& =& (q^{2m-1} + q^{-2m+1}) +(q^{2m-3} + q^{-2m+3}) 
+...+(q +q^{-1})\\
\nonumber
&=& z^{2m-1} - \binom{2m-2}{1}z^{2m-3}+ \binom{2m-3}{2}z^{2m-5}+...\\
&+& (-1)^r\binom{2m-r -1}{r}z^{2m-2r-1} +...+(-1)^{m-1}mz
\end{eqnarray}
Explicitly ( noting that in contrast with $(2.16)$ there is now $1$ 
on the right side below )
\begin{equation}
\Sigma _{N-2} =1 \qquad (N=3,5,...)
\end{equation}
  where
\begin{eqnarray}
\nonumber
   \Sigma _1 &=& z \\
\nonumber
  \Sigma _3 &=& z^3 - 2z \\
\nonumber
  \Sigma _5 &=& z^5 - 4z^3+ 3z \\
\nonumber
\Sigma _7& = &z^7 - 6z^5+ 10z^3 - 4z \\
  \Sigma _9& =& z^9 - 8z^7+21z^5 -20z^3 + 5z
\end{eqnarray}
and so on.

Comparing $(2.12)$ and $(2.13)$ one sees that $(B)$ is obtained from 
$(A)_1$, with now a nonzero
right side, as
\begin{equation}
  (B): \qquad S_n =2 \qquad (n= \frac{N}{2} = 2,3,...)
\end{equation}
Without trying to be exhaustive let us first point out some simple 
possibilities, particularly
when $(2.1)$ is obtained for $q$ a $\it  {root}$ $\it  {of}$ $\it 
{unity}$. ( See the relevant
remarks in the concluding  section.)

For $\hat{o}(3)$ one has from $(2.23)$
\begin{equation}
\Sigma _1 =q+q^{-1} =1,\qquad q=\frac{1}{2}+i\frac{\sqrt 3}{2}= 
e^{i\frac{\pi}{3}}, \qquad q^6 =1
\end{equation}

For $\hat{o}(4)$ one has from $(2.16)$ and $(2.20)$
\begin{equation}
Y = q^2 +q^{-2} =0,\qquad q=e^{i\frac{\pi}{4}}, \qquad q^8 =1
\end{equation}

Moreover, whenever $c_0 =0$ from $(2.19)$ i.e. for
\begin{equation}
  \hat {o}(4),\hat {o}(10),\hat {o}(12),..
\end{equation}
$Y$ can be factorized and $(2.27)$ is a solution.

For $\hat {o}(6)$ the two roots of
$$S_2 = (Y-1)(Y+2) =0$$
correspond to
\begin{equation}
  q^{12}=1, \qquad  q^{4}=1
\end{equation}

For $\hat {o}(8)$  the roots of $S_3$ give
\begin{equation}
  q^{6}=1, \qquad  q^{16}=1
\end{equation}

For $\hat {o}(10)$ and $\hat {o}(12)$ ( apart from $(2.27)$) one has 
to solve $\it {cubic}$ and
$\it {quartic}$ equations repectively. We do not present this standard algebra 
here. For odd $N$
one has again  a cubic in $z$ for $\hat {o}(5)$.

For $\hat {o}(14)$ onwards for $\hat {o}(2n)$ one has polynomials of 
sixth and higher degrees
(already for $Y$ before obtaining $q$ from $(2.15)$). Hence one needs 
$\it {hyperelliptic}$
functions for $Y$.

  For $\hat {o}(7)$ and $\hat {p}(10)$ one has quintics and  $\it 
{elliptic}$ solutions respectively
for $z = (q+q^{-1})$  and $ Y = (q^2+q^{-2})$. For higher dimensions one again
encounters hyperelliptic functions here.

  It is known $[6,7]$ that the general case on  a complex field
\begin{equation}
f(x)= a_0 x^n +a_1 x^{n-1} + ... + a_n = 0
\end{equation}
can be solved in terms of theta functions of zero arguments and the 
period matrix of the
hyperelliptic curves
\begin{equation}
  F^2 = x(x-1)f(x) , \qquad F^2 = x(x-1)(x-2)f(x)
\end{equation}
  for  odd $n$ and even $n$  rspectively.

   For quintics $[7]$ one can, alternatively, implement further 
successive changes of variables
  (Tschirnhausen transformations) to obtain standard forms (the 
Bring-Jerrard quintic or the
Brioschi quintic leading to the Jacobi sextic)  which can  be solved 
directly using elliptic
functions. All this is however very  complicated.

     The coeffficients $a_i$ of $(2.31)$ are very special ones 
(binomial integers) for our case.
What  special ( hopefully simplifying ) features might they induce in 
the corresponding elliptic
and hyperelliptic functions ? An answer to this question is beyond 
the  scope of this paper.

  Let us contemplate the simplest case, that of $\hat {o}(3)$ with 
$(2.26)$ giving

\begin {equation}
\hat R = I -\begin{vmatrix} 0 & 0 & 0 & 0 & 0 & 0 & 0 & 0 & 0 \\ 0 & 
0 & 0 & 0 & 0
& 0 & 0 & 0 & 0 \\ 0 & 0 & q^{-1} & 0 & q^{-\frac{1}{2}} & 0 & 1 & 0 
& 0 \\  0 & 0 & 0 & 0 & 0 & 0
& 0 & 0 & 0\\ 0 & 0 & q^{-\frac{1}{2}} & 0 & 1 & 0 & q^{\frac{1}{2}} 
& 0 & 0\\ 0 & 0 & 0 & 0 & 0 &
0 & 0 & 0 & 0\\ 0 & 0 & 1 & 0 & q^{\frac{1}{2}} & 0 & q & 0 & 0 \\ 0 
& 0 & 0 & 0 & 0 & 0 & 0 & 0 &
0 \\ 0 & 0 & 0 & 0 & 0 & 0 & 0 & 0 & 0 \end {vmatrix}
\end{equation}

where $ q=e^{i\frac{\pi}{3}}$.

Remarkably this satisfies

$${\hat R}_{12} {\hat R}_{23} {\hat R}_{12} =   {\hat R}_{23}{\hat 
R}_{12} {\hat R}_{23}$$

along with

$${\hat R}^2 =I$$

\subsection {A non-existence theorem:}

   For the standard cases ( see $(2.3),((2.4),(2.5)$ )

$$ (\hat R )_{q=1}=P, \qquad {\hat R}^2 =P^2 =I$$

  For the nonstandard Jordanian case ( see Refs.$8,9$ citing basic 
sources ) considering again
vector representations

\begin{equation}
\hat R = F^{-1}PF, \qquad {\hat R}^2 =I
\end{equation}
where $F$ is obtained through a "contraction" $[8,9]$. Thus the 
Yang-Baxter matrix
\begin{equation}
R= P\hat R \quad = (PF^{-1}P)F \quad  = {F}^{-1}_{21}F
\end{equation}
  is a "Drinfeld twist" of unity. This leads to various interesting 
features $[8,9]$ making
triangularity inherent without rendering $\hat R$ trivial. Can our 
constructions above ( for $\eta
=0$) be expressed as a conjugation of $P$ as in $(2.34)$ ? A priori 
such a possibility cannot be
discarded.

  However, using  our diagonalizers ( see App.$B$ of Ref.$1$ for 
explicit constructions ) one can
prove quite  simply and generally that $\it {no}$ invertible $F$ 
exists that can realize $(2.34)$.

     It is sufficient to to consider $\hat {o} (3)$. Higher dimensions 
can be treated in a strictly
parallel fashion. The essential result for us is that  the 
diagonalizer $M$ gives for $(1.26)$
\begin{equation}
M{\hat R}^{\pm 1}M^{-1}  = (-e^{\mp 2 \eta},1,,1,1,1,1,1,1,1)_{(diag)}
\end{equation}

For $\eta =0$, when ${\hat R}^2 =I$, one  thus obtains

\begin{equation}
M{\hat R}^{\pm 1}M^{-1}  = (-1,1,,1,1,1,1,1,1,1)_{(diag)} \equiv D
\end{equation}

  Now assume that an $F$ exists for our $\hat R$ satisfying $(2.34)$. Then

\begin{equation}
MF^{-1}PFM^{-1}  = D
\end{equation}
Defining
\begin{equation}
G=FM^{-1}
\end{equation}
since $D^2=I$
\begin{equation}
G=PGD
\end{equation}

The action of $P$ here for the $9\times 9$ case ( see ($2.3$) ) 
leaves the rows $(1,5,9)$ untouched
and interchanges the pairs of rows

                                $$ (2,4), (3,7), (6,8)$$

Now consider the action of $D$ after parametrizing $G$ in terms of 
parameters arbitrary to start
with. In rows $(1,5,9)$ the first element is constrained to be zero 
the others being unrestricted.
If row-$2$ is parametrized as

\begin{equation}
(a_1,a_2,a_3,a_4,a_5,a_6,a_7,a_8,a_9)
\end{equation}

then row-$4$ must be

\begin{equation}
(-a_1,a_2,a_3,a_4,a_5,a_6,a_7,a_8,a_9)
\end{equation}

Hence with arbitrary $a_1$, $r_i$ denoting the row-$i$,

\begin{equation}
(r_2 - r_4) = (2a_1,0,0,0,0,0,0,0,0)
\end{equation}

Similarly, in evident notations,
\begin{equation}
(r_3 - r_5) = (2b_1,0,0,0,0,0,0,0,0)
\end{equation}
\begin{equation}
(r_6 - r_8) = (2c_1,0,0,0,0,0,0,0,0)
\end{equation}

This evidently implies that the  determinant

\begin{equation}
\Delta G = 0
\end{equation}

Hence $G$ is not invertible. Hence neither is $F=GM$.

  This contradicts the assumption that an invertible $F$ exists giving 
$(2.34)$ for our $\hat R$.

\section {L- algebra (group-like central elements):}

  Before writing down the $RLL$- equations and the implied constraints 
we signal  the most
remarkable features to emerge in Secs.$3$ and $4$. We will study them 
mostly in the context of
the simplest case
$\hat {o}(3)$ i.e. $N=3$.

  ($1$): In the $L^+$ subalgebra one obtains $2N$ central, group-like 
elements constrained to
equality by the $RLL$-equations. There are $2N$ corresponding ones 
for the $L^-$ subalgebra.

($2$): In standard cases group-like elements are usually associated 
to "quantum determinants". But
our above-mentioned sets have no determinant-like structure at all. 
Each one is the sum of $N$
quadratic terms ( no negetive signs ).

($3)$: In the $3\times 3$ fundamental representation of the 
$L$-operators for $\hat {o}(3)$ these
elements are proportional to $I_3$. Consistently with their 
group-like property and centrality
they are proportional to $I_9$ for the $9\times 9$ coproduct 
representations. The explicit
verification of this involves remarkable cancellations. Iterated 
coproducts of course lead to
$I_{3^{2^{p}}}$ at the $p$-th stage.

($4$): In the standard cases we are used to the coproducts being 
reducible. Thus the $9\times 9$
coproducts  $\Delta L_{ij}^{\pm}$ for $SO_q (3)$ can be conjugated to 
block-diagonal forms
corresponding to the familiar irreducible components ($ 9\times 9 
\rightarrow 5\times  5 \oplus
3\times 3 \oplus 1\times 1$ )  or in terms of angular momenta ($ 
1\times 1 \rightarrow 2\oplus
1 \oplus 0$ ). But here one encounters obstructions in a systematic 
search for block-diagonalizations. This
is of course consistent with the central elements announced above 
being proportional to $I$. But since such
a search reveals special features of the generators this aspect of 
the $9\times 9$ coproduct mentioned
above will be treated explicitly in the next section.

Let us now formulate the $RLL$-constraints. The $FRT$-equations $[4]$ 
for the $L^{\pm}$ operators
are, in our notations,

\begin{equation}
\hat R L^{\pm}_2L^{\pm}_1 =L^{\pm}_2L^{\pm}_1\hat R
\end{equation}

\begin{equation}
\hat R L^{+}_2L^{-}_1 =L^{-}_2L^{+}_1\hat R, \qquad ({\hat R}^{-1} 
L^{-}_2L^{+}_1
=L^{+}_2L^{-}_1{\hat R}^{-1})
\end{equation}

Here $\hat R$ is a $N^2 \times N^2$ matrix satisfying $(1.21)$ and 
$L^{\pm}$ have each $N^2$
components $L^{\pm}_{ij},\quad (i,j = 1,2,..,N)$ arranged in a 
$N\times N$ matrix form with

\begin{equation}
L^{\pm}_2 =I_N \otimes L^{\pm}, \qquad L^{\pm}_1 = L^{\pm}\otimes I_N
\end{equation}

 From ($1.26$), with $P'_0$ given by ($1.17$),($1.18$),

\begin{equation}
{\hat R}^{±1} = I + {\lambda}_{\pm}P'_0 ,\qquad {\lambda}_{\pm} = - 
e^{\mp \eta}
\end{equation}

Here,from ($1.9$), ${\lambda}_{\pm}$ are the roots of

\begin{equation}
\lambda+{\lambda}^{-1} +([N\mp 1]\pm1) =0
\end{equation}
and , in particular, for $\hat o (3)$ of

\begin{equation}
\lambda+{\lambda}^{-1} +(q+1+q^{-1}) =0
\end{equation}

This simple change of notation $(\eta \rightarrow \lambda )$  will 
permit below a compact,
unified treatment of $L^+$ and $L^-$ due to the symmetry

\begin{equation}
(\lambda_+,L^+,L^-) \rightleftharpoons (\lambda_-,L^-,L^+)
\end{equation}

This is a special feature of our class. We will often suppress below 
the superscripts and
subscripts of $L$ and $\lambda$ respectively when dealing with 
($3.1$) and write $L_{ij}$ and
$\lambda$.

  From  ($3.1$) and  ($3.4$) one has

\begin{equation}
P'_0  L^{\epsilon}_2L^{\epsilon}_1 =L^{\epsilon}_2L^{\epsilon}_1P'_0, 
\qquad (\epsilon = \pm)
\end{equation}

In these equations $\lambda _{\pm}$ do $\it not$ appear explicitly. 
They do however appear for
($3.2$).

\subsection { $L^+$ and $ L^-$ subalgebras:}

  The group-like elements belong to these subalgebras since the 
coproducts are defined separately
for each. The $RLL$ constraints lead to ($2N^2 -1$) equations for 
each subalgebra which separate
into  three subsets of $N(N-1), N(N-1),(2N-1)$ respectively. The 
total number is easily understood
as follows.

  The diagonalizer $M$ ( see App.$B$, Ref.$1$ ) gives

\begin{equation}
M P'_0 M^{-1} \approx (1,0,0,....,0)_{(diag)}
\end{equation}

Conjugating the factors on both sides of ($3.8$) by $M$ only the 
row-$1$ survives on the left and
only the col.-$1$ on the right. These have each $N^2$ elements and 
one element in common. Hence the
result.

  It is convenient to start by deriving some results in a form valid, 
more generally, for the whole
algebra as follows. We consider $\hat o(N)$ for definiteness, the 
modifications for $\hat p (N)$
are evident.

  Along with ($1.17$) and ($3.3$) ,for $N \geq 3$, we write

\begin{equation}
L_2 L'_1 = \sum _{i,j,k,l} L_{kl}L'_{ij}(ij)\otimes (kl)
\end{equation}

where

$$ L_2L'_1 = L_2^{\epsilon}L_1^{{\epsilon}'}, \qquad (\epsilon, {\epsilon}') =
(+,+),(-,-),(+,-),(-,+)$$

One can show that (with $i' = N-i+1$)

\begin{equation}
P'_0 L_2 L'_1 = \sum _{i,,k,l} q^{-{\rho}_i}S^{(1)}_{lk}(ik)\otimes (i'l)
\end{equation}

\begin{equation}
  L'_2 L_1P'_0 = \sum _{i,,k,l} q^{-{\rho}_i}S^{(2)}_{lk}(ki)\otimes (li')
\end{equation}

where

\begin{equation}
S^{(1)}_{lk} = \sum _{j} q^{-\rho_{j}}L_{j'l}L'_{jk}
\end{equation}

\begin{equation}
S^{(2)}_{lk} = \sum _{j} q^{-\rho_{j}}L'_{lj'}L_{kj}
\end{equation}

We now go back to our subalgebras by setting $\epsilon = {\epsilon}'$ 
i.e. $L' =L$.

For $l\neq k'$ ( $l+k \neq N+1$ ) one obtains from ($3.8$),due to the 
structure of $P_0$,

\begin{equation}
S^{(1)}_{lk} = \sum _{j} q^{-\rho_{j}}L_{j'l}L_{jk} =0
\end{equation}
and

\begin{equation}
S^{(2)}_{lk} = \sum _{i} q^{-\rho_{i}}L_{li'}L_{ki} =0
\end{equation}

Each set ($S^{(1)},S^{(2)}$) corresponds to $N(N-1)$ equations and

\begin{equation}
S^{(1)} \rightarrow S^{(2)} \rightleftharpoons L_{ij} \rightarrow L_{ji}
\end{equation}

For
       $$l+k=N+1$$
one obtains ($i$ and $j$ assuming each value $\it independently $ 
with $i+i'=j+j' =N+1$ )
\begin{equation}
   q^{-\rho_i}S^{(1)}_{ii'} =  q^{-\rho_j}S^{(2)}_{jj'} \qquad 
(i=1,..N ;j=1,..N)
\end{equation}

The equality of these $2N$ quadratic expressions give ($2N-1$) 
equations. We will denote this set
as $\hat S_3$. This provides  the group-like central elements. We 
will study this set in detail
for
$\hat o (3)$. For  $N=3$ one obtains from ($3.18$)
\begin {equation}
q^{-\frac{1}{2}}S^{(1)}_{13} =
S^{(1)}_{22}=q^{\frac{1}{2}}S^{(1)}_{31}=q^{-\frac{1}{2}}S^{(2)}_{13}
=S^{(2)}_{22}=q^{\frac{1}{2}}S^{(2)}_{31}
\end{equation}

Or explicitly ( with $L$ all $L^+$ or all $L^-$ below )

\begin{eqnarray}
\nonumber
(\hat S_3)&:& L_{11}L_{33} +q^{-\frac{1}{2}}L_{21}L_{23}+q^{-1}L_{31}L_{13} \\
\nonumber
&=&q^{\frac{1}{2}}L_{12}L_{32} +L_{22}L_{22}+q^{-\frac{1}{2}}L_{32}L_{12}\\
\nonumber
&=&qL_{13}L_{31} +q^{\frac{1}{2}}L_{23}L_{21}+L_{33}L_{11} \\
\nonumber
&=&L_{11}L_{33}+q^{-\frac{1}{2}}L_{12}L_{32}+ q^{-1}L_{13}L_{31} \\
\nonumber
&=&q^{\frac{1}{2}}L_{21}L_{23} +L_{22}L_{22}+q^{-\frac{1}{2}}L_{23}L_{21} \\
&=&qL_{31}L_{13} +q^{\frac{1}{2}}L_{32}L_{12}+ L_{33}L_{11}
\end{eqnarray}

We denote the sets ($3.15$) and ($3.16$) by $\hat S_1$ and $\hat S_2$ 
respectively. For $N=3$
the first set is

$$\hat S_1:$$

$$ q^{-\frac{1}{2}}L_{31}L_{11}+L_{21}L_{21}+q^{\frac{1}{2}}L_{11}L_{31} =0$$

$$ q^{-\frac{1}{2}}L_{32}L_{11}+L_{22}L_{21}+q^{\frac{1}{2}}L_{12}L_{31} =0$$

$$ q^{-\frac{1}{2}}L_{31}L_{12}+L_{21}L_{22}+q^{\frac{1}{2}}L_{11}L_{32} =0$$

$$ q^{-\frac{1}{2}}L_{33}L_{12}+L_{23}L_{22}+q^{\frac{1}{2}}L_{13}L_{32} =0$$

$$ q^{-\frac{1}{2}}L_{32}L_{13}+L_{22}L_{23}+q^{\frac{1}{2}}L_{12}L_{33} =0$$
\begin{equation}
  q^{-\frac{1}{2}}L_{33}L_{13}+L_{23}L_{23}+q^{\frac{1}{2}}L_{13}L_{33} =0
\end{equation}

   The set $\hat S_2$ is obtained immediately from $\hat S_1$ via ($3.17$).

We now show how the sets ($\hat S_1,\hat S_2,\hat S_3$) imply the two 
basic properties of the
members of $\hat S_3$, that they  are ($1$) cnetral and ($2$) group-like.

($1$): Exploiting systematically the sets $\hat S_i$ one can pass 
through different chains of
intermediate steps. One possible sequence is as follows.

\begin{equation}
q^{\frac{1}{2}}L_{11}L_{23}L_{21} = q^{-\frac{1}{2}}L_{12}L_{32}L_{11}
-q^{-\frac{1}{2}}L_{13}L_{21}L_{21} + q^{\frac{1}{2}}L_{12}L_{12}L_{31}
\end{equation}

\begin{equation}
qL_{11}L_{13}L_{31} = q^{-1}L_{13}L_{31}L_{11}
+q^{-\frac{1}{2}}L_{13}L_{21}L_{21} - q^{\frac{1}{2}}L_{12}L_{12}L_{31}
\end{equation}

Summing these with $L_{11}L_{33}L_{11}$ one obtains ( using again 
($3.20$) in the last step )
\begin{eqnarray}
\nonumber
L_{11} ( qL_{13}L_{31} +q^{\frac{1}{2}}L_{23}L_{21}+L_{33}L_{11} )
\nonumber
&=&(L_{11}L_{33}+q^{-\frac{1}{2}}L_{12}L_{32}+ q^{-1}L_{13}L_{31}) L_{11} \\
&=&( qL_{13}L_{31} +q^{\frac{1}{2}}L_{23}L_{21}+L_{33}L_{11} )L_{11}
\end{eqnarray}

Thus $L_{11}$, and similarly each $L_{ij}$ can be shown to commute 
with the members of $\hat
S_3$. Hence the latter are central.

($2$): The rule for coproducts $[4]$ ( with $L$ either $L^+$ or $L^-$ 
throughout ) is

\begin{equation}
\Delta L_{ij} = \sum _{k}L_{ik}\otimes L_{kj}
\end{equation}

Hence
\begin{equation}
\Delta L_{ab}\Delta L_{cd} = \sum _{i,j}L_{ai}L_{cj}\otimes L_{ib}L_{jd}
\end{equation}

Now let us start with the first member of $\hat S_3$ and compute the sum

\begin{equation}
\Sigma = \Delta L_{11}\Delta L_{33} + q^{-\frac{1}{2}}\Delta 
L_{21}\Delta L_{23}+
q^{-1}\Delta L_{31}\Delta L_{13}
\end{equation}

Collecting terms and systematically implementing $\hat S_1,\hat S_2$ 
many terms cancel leaving

\begin{equation}
\Sigma = ( L_{11} L_{33} + q^{-\frac{1}{2}} L_{21} L_{23}+q^{-1} L_{31} L_{13})
\otimes ( L_{11} L_{33} + q^{-\frac{1}{2}} L_{21} L_{23}+q^{-1} L_{31} L_{13})
\end{equation}

Thus this and similarly the other members of $\hat S_3$ are group-like.

For $\hat o(4)$ the $8$ members of $\hat S_3$ satisfy

$$ q^{-1}S^{(1)}_{14} = S^{(1)}_{23} = S^{(1)}_{32} = qS^{(1)}_{41}$$

\begin{equation}
= q^{-1}S^{(2)}_{14} = S^{(2)}_{23} = S^{(2)}_{32} = qS^{(2)}_{41}
\end{equation}

Their centrality and group -like property  can be established in 
strict analogy to the $\hat
o(3)$ case. Moreover they indicate how the generalization for higher 
$N$ can be carried out. Our
presentation here will be limited to $\hat o(3)$.

\subsection{Beyond the $L^{\pm}$ subalgebras:}

  We now consider the "mixed" case ($3.2$). The major new feature now 
is the explicit involvement
of $\lambda _{\pm}$ (see ($3.4$),($3.5$)) in the constraints

\begin{equation}
(I + \lambda_{+}P'_{0})L_2^+L_1^- =L_2^-L_1^+(I + \lambda_{+}P'_{0})
\end{equation}

\begin{equation}
(I + \lambda_{-}P'_{0})L_2^-L_1^+ =L_2^+L_1^-(I + \lambda_{-}P'_{0})
\end{equation}

But even when $P'_{0}$ does not contribute, namely at

$$ (ij)\otimes (kl), \qquad (kl) \neq (i'j')$$

one obtains simple  but probing constraints. Retaining only such rows 
and columns one obtains a
"reduced" matrix of $N(N-1)\times N(N-1)$ dimensions. For $N=3$ this 
corresponds to the
suppression of rows and columns ($3,5,7$) leaving a $6\times  6$ 
matrix. Using for this reduced
case, for all $N$, the subscript $r$ one extracts from ($3.30$),($3.31$)

\begin{equation}
(L^{\epsilon}_2L^{\epsilon'}_1)_{(r)}=(L^ {\epsilon'}_2L^{\epsilon}_1)_{(r)}
\end{equation}

For $\epsilon = \epsilon'$ this is trivial. But not now and one can 
go further as follows.

Since $\lambda$ satisfies a quadratic equation one can linearize all 
polynomials in $\lambda$
using $\lambda_+$ and $\lambda_- = (\lambda_+)^{-1}$. The symmetry of 
($3.30$),($3.31$),

\begin{equation}
(\lambda_{+} \rightarrow \lambda_{-}) \rightleftharpoons ( L^+ 
\rightarrow L^- )
\end{equation}

indicates the parametrization where the $\lambda$-dependence is 
explicitly ( and only ) in the
coefficient as
\begin{equation}
L^{\epsilon}_{ij} = A_{ij}+ \lambda_{\epsilon} B_{ij}
\end{equation}

  Now injecting ($3.34$) in ($3.32$) oe obtains for each element of of 
the {\it reduced matrix}

( $L^{\pm}_{ab},L^{\pm}_{cd}$ and so on )

\begin{equation}
(\lambda_{+} - \lambda_{-}) ( A_{ab}B_{cd} -  B_{ab}A_{cd} ) = 0, 
\qquad (\lambda_{+} \neq
\lambda_{-})
\end{equation}

Thus the $\lambda$-dependence is simply factored out for this redued 
martix. Even when $L^+$ and
$L^-$ are considered in the  context of the respective subalgebras 
they must satisfy ($3.35$).
This will indeed be found to be  the case in the  explicit 
realizations of the following
section.

  We now come to parts where both $I$ and $P'_0$ contribute and hence 
$\lambda$ is directly
involved. Instead of $\hat S_3$ of ($3.20$) one now has for $\hat 
o(3)$ $9$ relations of the type

\begin{eqnarray}
\nonumber
\lambda_{+} \biggl(( qL'_{31}L_{13} +q^{\frac{1}{2}}L'_{32}L_{12}+ 
L'_{33}L_{11} ) &-&(
qL_{13}L'_{31} +q^{\frac{1}{2}}L_{23}L'_{21}+ L_{33}L'_{11} )\biggr) \\
  &=& q( L_{33}L'_{11} - L'_{33}L_{11})
\end{eqnarray}

   Instead of $\hat S_1,\hat S_2$ one now has equations of the type

$$ \lambda_+ 
(q^{-\frac{1}{2}}L_{31}L'_{11}+L_{21}'L_{21}+q^{\frac{1}{2}}L_{11}L'_{31}) 
$$

\begin{equation}
= q^{\frac{1}{2}}(L'_{31}L_{11} - L_{31}L'_{11}) = (L'_{21}L_{21} - 
L_{21}L'_{21}) =
q^{-\frac{1}{2}}(L'_{11}L_{31} - L_{11}L'_{31})
\end{equation}

For  $\epsilon = \epsilon'$ one recovers the results for  the 
subalgebras. We have obtained a
systematic formulation of the full set of $81$ constraints for $\hat 
o(3)$ exploiting certain
symmetries. This will not be reproduced here. The generalizations of 
the results of this
subsection for $N>3$ can be obtained fairly systematically.

\subsection { From $L^{\pm}$ to $L(\theta)$:}

  Since our $\hat R$ satisfies a quadratic equation ($1.24$) all the 
three FRT equations
(($3.1$),($3.2$)) can be  condensed into a {\it single } one by 
defining in analogy to

\begin{equation}
\hat R(\theta) = \frac {e^{\eta+\theta}\hat R - e^{-\eta-\theta}{\hat
R}^{-1}}{e^{\eta+\theta}-e^{-\eta-\theta}}
\end{equation}

\begin{equation}
L(\theta) = \frac {e^{\eta+\theta}L^+ - 
e^{-\eta-\theta}L^-}{e^{\eta+\theta}-e^{-\eta-\theta}}
\end{equation}

It can be shown that $[10,11]$

\begin{equation}
\hat R(\theta - \theta') L_{2}(\theta) L_{1}(\theta')= L_{2}(\theta') 
L_{1}(\theta)\hat R(\theta
-\theta')
\end{equation}

contains effectively all the three FRT equations. One can write 
($3.38$) i.e. ($1.10$) as

\begin{equation}
\hat R(\theta) = I - \frac {sinh \theta}{sinh (\eta +\theta)} P'_0
\end{equation}

\section {Fundamental and coproduct representations :}

    Here we will study the $\hat o (3)$ fundamental ($3\times 3$) and 
the coproduct ($9\times 9$)
representations. They illustrate the significance of the remarks 
($3,4$) at the beginning of
Sec.$3$. This will be commented upon at the end. Specific symmetries 
of the matrices obtained
will be displyed. They might be helpful in a more systematic study of 
representations.

  The general prescription for the fundamental representations 
($N\times N$ blocks for all $N$)
is as follows. (See App.$B$ of Ref.$11$ for a systematic presentation 
citing basic  sources.)

\begin{equation}
L(\theta) = \hat R(\theta) P, \qquad ( L(\theta))_{q\rightarrow 
\pm{\infty}} = L^{\pm} =  \hat
R^{\pm 1} P
\end{equation}

  A Hopf algebra can be defined [$4,10$] for $L$ using

\begin{equation}
(L^{\pm})^{-1} = PL^{\mp}P = (L^{\mp})_{21}
\end{equation}

 From ($3.4,3.6$) and ($4.1$) one obtains ( for the fund. repr. )

$$ L = \begin{vmatrix} L_{11} &L_{12} &L_{13} \\ L_{21} &L_{22} 
&L_{23} \\ L_{31} &L_{32}
&L_{33} \end{vmatrix} $$

$$ L_{11} = \begin{vmatrix} 1 &0 & 0 \\ 0 &0 &0 \\0 &0&\lambda\end{vmatrix}, 
\qquad
L_{12}=\begin{vmatrix} 0 &0 & 0 \\ 1 &0 &0 \\0 
&q^{-\frac{1}{2}}\lambda&0\end{vmatrix}, \qquad
L_{13}=\begin{vmatrix} 0 &0 & 0 \\ 0 &0 &0 
\\(1+q^{-1}\lambda)&0&0\end{vmatrix}$$

$$ L_{21} = \begin{vmatrix} 0 &1 & 0 \\ 0 &0 &q^{\frac{1}{2}}\lambda \\0
&0&0\end{vmatrix}, \qquad L_{22}=\begin{vmatrix} 0 &0 & 0 \\ 0 
&(1+\lambda) &0 \\0
&0&0\end{vmatrix}, \qquad L_{23}=\begin{vmatrix} 0 &0 & 0 \\ 
q^{-\frac{1}{2}}\lambda &0 &0
\\0&1&0\end{vmatrix}$$

\begin{equation}
  L_{31} = \begin{vmatrix} 0 &0 & (1+q \lambda) \\ 0 &0 &0 \\0 
&0&0\end{vmatrix}, \qquad
L_{32}=\begin{vmatrix} 0 &q^{\frac{1}{2}}\lambda & 0 \\ 0 &0 &1 \\0
&0&0\end{vmatrix}, \qquad L_{33}=\begin{vmatrix} \lambda &0 & 0 \\ 0 &0 &0
\\0&0&1\end{vmatrix}
\end{equation}

One sets $\lambda = (\lambda _+,\lambda_-,\lambda (\theta))$ for $ L= 
(L^+,L^-,L(\theta))$
respectively. One obtains $\lambda (\theta)$ directly from ($3.39$).

Note that ($3.34$) is evidently satisfied with

$$ L^{\pm}_{11} = \begin{vmatrix} 1 &0 & 0 \\ 0 &0 &0 \\0
&0&0\end{vmatrix}+\lambda_{\pm}\begin{vmatrix} 0 &0 & 0 \\ 0 &0 &0 
\\0 &0&1\end{vmatrix} \equiv
A_{11} +\lambda_{\pm} B_{11} $$
and so on.

   Now we consider the $9\times 9$ coproducts of ($4.3$) given by

     $$\Delta L_{ij}^{\pm} = \sum_{k}L_{ik}^{\pm}\otimes L_{kj}^{\pm}$$

We will not present the easily obtained $9\times 9$ matrices but the 
symmetries they exhibit for
the reason mentioned before. Define

($r_1,r_2$) : reflections about the diagonal aand the antidiagonal respectively

$f:  (q\rightarrow q^{-1}) (r_2r_1), \qquad f(A) \equiv fA$

Then in terms of $3\times 3$ blocks $A_{ij},...,E_3$ (not  exhibited 
here) one obtains with ($ij
= 11,12,13$) and ($i'j'= 33,32,31$) respecttively ( and noting that 
$\lambda$ is invariant for
$q\rightarrow q^{-1}$ )

\begin{equation}
\Delta L_{ij} =\begin{vmatrix} A_{ij} &0 &0 \\ B_{ij} &0 &0 \\ C_{ij} 
&q^{-\frac{1}{2}}\lambda
B_{ij} &\lambda A_{ij} \end{vmatrix}, \qquad f\Delta L_{ij} = \Delta L_{i'j'}
=\begin{vmatrix}\lambda f A_{ij} &q^{\frac{1}{2}}\lambda fB_{ij} &fC_{ij}
\\ 0 &0& fB_{ij} \\0 &0 & fA_{ij} \end{vmatrix}
\end{equation}

\begin{equation}
\Delta L_{21} =\begin{vmatrix} 0 &D_1 &0 \\ q^{-\frac{1}{2}}\lambda D_2 &D_3
&q^{\frac{1}{2}}\lambda D_1 \\ 0 &D_2 &0\end{vmatrix}, \qquad f\Delta 
L_{21} = \Delta L_{23}
=\begin{vmatrix} 0 &fD_2 &0 \\ q^{-\frac{1}{2}}\lambda fD_1 &fD_3
&q^{\frac{1}{2}}\lambda fD_2 \\ 0 &fD_1 &0\end{vmatrix}
\end{equation}

\begin{equation}
\Delta L_{22} =f\Delta L_{22} =\begin{vmatrix} 0 &E_1 &0 \\ 
q^{-\frac{1}{2}}\lambda fE_1 &E_3
&q^{\frac{1}{2}}\lambda E_1 \\ 0 &fE_1 &0\end{vmatrix}
\end{equation}

   Having displayed the symmetries we now study in more detail the 
three generators $\Delta
L_{ii}$. For the standard cases the $L_{ii}$ can be obtained directly 
in diagonal forms for
irreducible representations [$4$] and through appropriate 
conjugations for reducible ones. For
our $3\times 3$ representation also they are diagonal with their sum 
proportional to $I$. But for the
$9\times 9$ coproducts above they  are not diagonal. They do not 
commute mutually and hence
cannot be diagonalized simultaneously. To better understand the 
structure encountered let us
try to diagonalize the sum ($\sum_i \Delta L_{ii}$) . Define

\begin{equation}
\mu = (q^{\frac{1}{2}}-q^{-\frac{1}{2}}), \quad z= 
(q^{\frac{1}{2}}+q^{-\frac{1}{2}}), \quad k =
(q+4+q^{-1})^{-\frac{1}{2}}
\end{equation}

Set

\begin{equation}
\sqrt{2}N = \begin{vmatrix} 0 & 1 & 0 & 1 & 0 & 0 & 0 & 0 & 0 \\ 0 & 
1 & 0 & -1 & 0
& 0 & 0 & 0 & 0 \\ 0 & 0 & 0 & 0 & 0 & 1 & 0 & -1 & 0 \\  0 & 0 & 0 & 0 & 0 & 1
& 0 & 1 & 0\\ 0 & 0 & 1 & 0 & 0 & 0 &-1 & 0 & 0\\ 0 & 0 &zk & 0 & 2k
& 0 & zk & 0 & 0\\ 0 & 0 & \sqrt{2}k & 0 &-\sqrt{2}zk & 0 &\sqrt{2}k 
& 0 & 0 \\ \sqrt{2} & 0 & 0
& 0 & 0 & 0 & 0 & 0 & 0 \\ 0 & 0 & 0 & 0 & 0 & 0 & 0 & 0 & 
\sqrt{2}\end {vmatrix}
\end{equation}

This is orthogonal i.e.

               $$ N^{-1} =N^T$$

One obtains
\begin{equation}
N (\Delta L_{11}+\Delta L_{22}+\Delta L_{33})N^{-1} = \lambda 
(-y,y,y,-y,3,3,-y,-y,-y)_{(diag)}
,\qquad y= (q+1+q^{-1})
\end{equation}

  The eigenvalues can be permuted through evident supplementary conjugations.

  We have thus diagonalized the sum. Now let us look at the component 
terms. One has ( denoting by
$(bd)$ a block diagonal structure )

\begin{equation}
N (\Delta L_{ii})N^{-1} = 
(\alpha_{i},\beta_{i},\gamma_{i},\delta_{i})_{(bd)} \qquad (i =1,2,3)
\end{equation}

where
\begin{equation}
\alpha_{1} =\frac{1}{2} \begin{vmatrix} 
1&1\\-1&-1\end{vmatrix},\qquad \alpha_{3}
=\frac{\lambda^{2}}{2}\begin{vmatrix} 1& -1\\1&-1\end{vmatrix}, 
\qquad \alpha_{2}
=\frac{1}{2}\begin{vmatrix} \lambda^{2}+1&\lambda^{2}
-1\\-\lambda^{2}+1&-\lambda^{2}-1\end{vmatrix}
\end{equation}

\begin{equation}
\beta_{1} =\frac{\lambda^{2}}{2} \begin{vmatrix} 
-1&-1\\1&1\end{vmatrix},\qquad \beta_{3}
=\frac{1}{2}\begin{vmatrix} -1& 1\\-1&1\end{vmatrix}, \qquad \beta_{2}
=\frac{1}{2}\begin{vmatrix}-\lambda^{2}-1&\lambda^{2}
-1\\-\lambda^{2}+1&\lambda^{2}+1\end{vmatrix}
\end{equation}

\begin{equation}
\delta_{1} = \begin{vmatrix} 1&0\\0&\lambda^{2}\end{vmatrix},\qquad \delta_{3}
=\begin{vmatrix}\lambda^{2} & 0\\0&1\end{vmatrix}, \qquad \delta_{2}
=\begin{vmatrix}0&0\\0&0\end{vmatrix}
\end{equation}

\begin{equation}
\gamma _{1} =\frac{\lambda}{2}\begin{vmatrix} 3 &\mu k &\sqrt{2}((2+q^{-1})k \\
\mu k &3 z^2 k^2 &\sqrt{2}((2+q^{-1})\mu k^2 \\ 
-\sqrt{2}((2+q)k&-\sqrt{2}((2+q)\mu k^2
&-2(z^2 -1)k^2\end{vmatrix}
\qquad
\end{equation}

\begin{equation}
\gamma _{3} =\frac{\lambda}{2}\begin{vmatrix} 3 &\mu k &-\sqrt{2}((2+q)k \\
\mu k &3(1- 2 k^2) &-\sqrt{2}((2+q)\mu k^2 \\ 
\sqrt{2}((2+q^{-1})k&\sqrt{2}((2+q^{-1})\mu k^2
&-2(1-3k^2)\end{vmatrix}
\end{equation}

\begin{equation}
\gamma _{2} =\frac{\lambda}{2}\begin{vmatrix} 0 &-2\mu k &\sqrt{2}\mu zk \\
-2\mu k &12 k^2 &\sqrt{2}\mu^2 z k^2 \\ \sqrt{2}\mu zk&\sqrt{2}\mu^2 z k^2
&-2(z^2 -1)z^2k^2\end{vmatrix}
\end{equation}

  Note that

\begin{equation}
\alpha _1^2 = \quad \alpha _3^2 = \quad \beta _1^2 = \quad \beta _3^2 = \quad
\begin{vmatrix}0&0\\0&0\end{vmatrix}
\end{equation}

  Such nilpotent matrices are {\it nondiagonalizable}.

It has been explicitly verified that not only $N (\Delta L_{ij}) 
N^{-1}$ $(i\neq j)$ are not
correspondingly block diagonalized but {\it all} their nonzero 
elements lie systematically
outside the blocks arising for $N (\Delta L_{ii}) N^{-1}$. One can 
examine larger blocks, say,
$(6\otimes 6 \oplus 3\otimes 3 )$ after permuting the $\gamma _i$ and 
$\delta_i$ blocks. But one
finds that the whole $9 \times 9$ space is needed for $N (\Delta 
L_{ij}) N^{-1}$. These results
will not be displayed here though particularly for $q=1$ they aquire 
relatively simple forms.
One can implement further conjugations and permutations but the 
essential features persist.

  For the $3\times 3$ representations all the members of $\hat S_3$ in 
$(3.20)$ say, for example,
\begin{equation}
  L_{11}L_{33} +q^{-\frac{1}{2}}L_{21}L_{23}+q^{-1}L_{31}L_{13} = \lambda I_3
\end{equation}

Also
\begin{equation}
  L_{11}+ L_{22}+L_{33} = (1+\lambda)I_3
\end{equation}

The members of $\hat S_3$ being group-like ($4.18$) gives for the 
$9\times 9$ coproducts
$\lambda ^2 I_9$. This has been verified explicitly. But ($ \sum_{i} 
\Delta L_{ii}$) behaves
quite differently as shown above.

  The obstructions encountered in reduction of the $9\times 9$ 
coproducts to smaller dimensional
irreducible components ( via block diagonalization in a fashion 
analogous, say, to the case of
$SO_q(3)$ ) is consistent with the central $ \hat S_3$ operators 
being proportional to $I$. But
our study of $\Delta L_{ii}$ reveals specific properties of these 
generators for higher
dimensional representations ( such as symmetries and 
nondiagonalizable blocks ). This can be
helpful in a more systematic study of representations. The symmetries 
displayed in
($4.4,4.5,4.6$) stem from  those of $P_0$ and hence should be 
significant more generally.

\section {Link invariants ( Turaev construction):}

\subsection {Construction of "enhanced" operators:}

  Given a matrix satisfying the braid equation the Turaev 
construction[$12$] of an enhanced Yang
Baxter operator ($EYB$) leads to explicit construction of invariants 
( invariant under Markov
moves of first and second types ) for oriented links. Such an 
enhanced system [$12,13$] consists
of a $N^2\times N^2$ braid matrix $\hat R$, an $N\times N$ matrix $f$ 
and elements( $a,b$), all
invertible, satisfying the relations

\begin{equation}
{\hat R}^{\pm 1}f\otimes f = f\otimes f {\hat R}^{\pm 1}
\end{equation}

\begin{equation}
tr_2({\hat R}^{\pm 1}f\otimes f) = a^{\pm 1}bf
\end{equation}

where one defines

\begin{equation}
tr_2 \biggl( \sum_{ijkl}c_{ij,kl}(ij)\otimes (kl) \biggr) = \sum_{ij} 
\biggl( \sum_k
c_{ij,kk}\biggr)(ij)
\end{equation}

  Let us first obtain ($f,a,b$) for our class of $\hat R$. Our 
spectral resolution and the
properties of the projector $P_0$ ( and hence of $P'_0$ defined in 
($1.17,1.18$) ) render the
constructions particularly transparent.

Define the ${\it diagonal}$ $N\times N$ matrices $f$  for the cases

$$ (1): \quad \hat o(N),\qquad N=2n+1, \quad n= 1,2,..$$
$$ (2): \quad \hat o(N),\qquad N=2n, \quad n= 2,3,..$$
$$ (3):\quad \hat p(N),\qquad N=2n, \quad n= 2,3,.. $$

respectively as follows

\begin{equation}
(1): \quad f = (q^{-(2n-1)}, 
q^{-(2n-3)},...,q^{-1},1,q,...,q^{(2n-3)},q^{(2n-1)})_{(diag)}
\end{equation}

\begin{equation}
(2): \quad f = (q^{-(2n-2)}, 
q^{-(2n-4)},...,q^{-2},1,1,q^2,...,q^{(2n-4)},q^{(2n-2)})_{(diag)}
\end{equation}

\begin{equation}
(3): \quad f = (q^{-2n}, 
q^{-(2n-2)},...,q^{-2},q^2,...,q^{(2n-2)},q^{2n})_{(diag)}
\end{equation}

   Note the following facts:

  ($1$): The $N$ diagonal elements of $f$, in each case, are the 
nonzero diagonal elements of
the corresponding $P'_0$ ( related to the projector as $P_0 = (tr 
P'_0)^{-1} P'_0$ ),the
remaining $N(N-1)$  diagonal elements of $P'_0$ being zero. Hence ( 
with upper and lower signs
for $\hat o(N)$ and $\hat p(N)$ respectively ),

\begin{equation}
T \equiv tr f \quad  = tr P'_0 \quad = [N \mp 1] \pm 1 \quad = 
-(\lambda _+ + \lambda_-) \quad
=e^{-\eta} +e^{\eta}
\end{equation}

where $e^{\pm \eta}$ are the  roots of

\begin{equation}
e^{2\eta} - T e^{\eta} + 1 =0
\end{equation}

($2$): The $N^2\times N^2$ matrix $f\otimes f$, in each case, has $1$ 
on rows {\it and} columns

$$ ( N,2N-1,3N-2,...,(N^2 -N +1) )$$

  These are precisely ones on  which $P'_0$ has nonzero elements. 
Hence directly ( without
further computations ) we obtain

\begin{equation}
  P'_0f\otimes f = f\otimes f  P'_0
\end{equation}

($3$): Using ($1.17,1.18, 5.3$) one obtains

\begin{equation}
tr_{2}( P'_0f\otimes f ) = f
\end{equation}

  Now from ($5.9,5.10$) it follows immediately

\begin{equation}
\hat R^{\pm 1}f\otimes f = (I - e^{\mp\eta} P'_0)f\otimes f = 
f\otimes f ( I - e^{\mp\eta} P'_0)
= f\otimes f\hat R^{\pm 1}
\end{equation}

and

\begin{equation}
tr_{2}(\hat R^{\pm 1}f\otimes f) = tr_{2}((I -  e^{\mp\eta} 
P'_0)f\otimes f) = (T - e^{\mp \eta})f
= e^{\pm \eta}f
\end{equation}

Hence

\begin{equation}
tr_{2}(\hat R^{\pm 1}f\otimes f) =  a^{\pm1}bf, \qquad (a= e^{\eta}, b=1)
\end{equation}

  Thus we have  obtained for our $\hat R$, in terms of $f$ introduced 
above,  the enhanced
operator

                       $$ ( \hat R,f, e^{\eta},1 )$$

  Our $f$ is srticly analogous  to those of Turaev for 
$SO_q(2n+1),SO_q(2n),Sp_q(2n)$
respectively. But whereas for the  standard cases $a$  ( $\alpha$ in 
the notation of [$12$] ) is
also a simple power of $q$, for us it involves the squareroot of a 
Laurent polynomial in $q$.
One obtains with $\delta = (1,2,0)$ for $\hat o (2n+1),\hat o 
(2n),\hat p (2n)$ respectively

\begin{equation}
T = (q^{-2n+\delta} + q^{-2n -2+\delta}+ ...+q^{2n+2-\delta} + 
q^{2n-\delta}+\delta )
\end{equation}
and
\begin{equation}
a^{\pm 1}= e^{\pm \eta}= \frac{1}{2}(T \pm \sqrt {T^2 -4})
\end{equation}

This is the crucial new aspect for our class of $\hat R$. For the 
simplest case $\hat o (3)$ one
obtains

\begin{equation}
a^{\pm 1}=\frac{1}{2} (q+1+q^{-1}) \pm \frac{1}{2} \biggl (
(q+3+q^{-1})(q-1+q^{-1})\biggr)^{\frac{1}{2} }
\end{equation}

giving for $q=1$,

\begin{equation}
a^{\pm 1}=\frac{3}{2}  \pm \frac{1}{2} \sqrt 5
\end{equation}

In general, for $q=1$, as for the standard case $f$ reduces to the 
$N\times N$ unit matrix but
as emphasized before our $\hat R$ remains nontrivial and

\begin{equation}
a^{\pm 1}=\frac{1}{2}N  \pm \frac{1}{2} \sqrt {N^2 - 4}, \qquad ( N= 3,4,...)
\end{equation}

  The discussion of Sec.$2$ shows that for solutions of ($2.10$)

  \begin{equation}
  \eta =0, \qquad    a=1
\end{equation}

implying a complex root of unity $q$ for $N=3$ but finally elliptic 
aand hyperelliptic ones as
$N$ increases. ( Overcrossings and undercrossings degenerate for 
$\hat R = {\hat R}^{-1}$.)

  [ Comparison of notations: The present author is often confused by 
different significances of
the same symbol ( and vice versa ) encountered elsewhere. The 
following points might be helpful
in our context.

        Turaev's $R$ [$12$] satisfying the braid equation ( his eqn. $1$ )

  $$ R_1R_2R_1 = R_2R_1R_2$$

  is our $\hat R$. Our $R$ is $P \hat R$ where

$$ P(ij)\otimes (kl) = (kj)\otimes (il), \quad  (ij)\otimes (kl)P = 
(il)\otimes (kj), \quad
P(ij)\otimes (kl)P = (kl)\otimes (ij)$$

and $R$ satisfies the Yang-Baxter equation

   $$ R_{12}R_{13}R_{23} = R_{23}R_{13}R_{12}$$

  It must also be clearly be noted that if $\sigma$ ( or $\tau$ ) is defined as

$$\sigma(ij)\otimes (kl) = \quad (kl)\otimes (ij) = \quad P(ij)\otimes (kl)P$$

then

$$ \sigma R = PRP = R_{21} = \hat R P$$

  does not stisfy the braid equation above ( satisfied by $\hat R$ ).

  Moreover, since

$$P(f\otimes f) P = f\otimes f , \qquad P^2 = I$$

  the condition ($5.1$) implies also ( for the YB-matrix $R$ )

$$R^{\pm 1}f\otimes f = f \otimes f R^{\pm 1}$$

  This form is presented in Sec.$15.2.2$ of Ref.$13$.But  $\hat R$ 
cannot be replaced by $R$ in
($5.2$). The symbol $I$ denotes $\sigma R$ in Sec.15.2.2 of Ref.$13$ 
and $-\sigma R$  at the end
of Sec.15.2.5. ]

\subsection{Link Invariants and skein relation:}

   We follow the presentation of Ref.$12$ and Sec.$15$ of Ref.$13$ 
with some changes of
notations. Let $\rho (\beta)$ be the representation of the braid 
$\beta$ associated to $\hat
R$ and let $\alpha (\beta)$ be the " augmentation homomorphism" 
changing by $\pm 1$ corresponding
to the actions of $T^{\pm1}$, the generators of the braid group.

  Define for our case ( with $b=1$ )

\begin{equation}
  \textsf{P} (\beta) = a^{-\alpha(\beta)} tr \bigl( \rho 
_{m}(\beta)\cdot f^{\otimes m}\bigr )
\end{equation}

$\rho_{m}$ being the endomorphism of $V^{\otimes m}$ associated to $\hat R$.

     Using appropriately the properties ($5.1$) and ($5.2$) of $f$ 
such a $\textsf{P} (\beta)$
can be shown to be Markov invariant and provide an invariant of 
oriented links. Markov moves
are defined, for example, in Sec.$15.1$ of Ref.$13$ and the proof of 
invariance of $ \textsf{P}
(\beta)$ is given in Sec.$15.2$ following  Ref.$12$.

  For an "unknot" ( no crossing) one has

\begin{equation}
  \textsf{P} (\bigcirc) = trf = T
\end{equation}

    Using standard notations ($\bar{L}_+,\bar{L}_-,\bar{L}_0$) 
corresponding to one point of the
projection of a braid differing by an overcrossing, undercrossing and 
nocrossing respectively
one obtains in our case ( following the steps below eqn.($6$), 
Sec.$15.2$ of Ref.$13$ )

$$x  \textsf{P} (\bar{L}_+) +y  \textsf{P} (\bar{L}_-)+ z  \textsf{P} 
(\bar{L}_0)$$
\begin{equation}
  = a^{-\alpha
(\beta)} tr \biggl( \bigl (x\hat R^{-1}+y \hat R +z I\bigr)\otimes 
id^{\otimes 2}\otimes
id^{\otimes (m- 2)}\cdot \rho_m (\beta)\cdot f^{\otimes m} \biggr)
\end{equation}

Now for our case
\begin{equation}
  e^{\eta} \hat R -  e^{-\eta} {\hat R}^{-1} = (  e^{\eta} - e^{-\eta} )I
\end{equation}

Hence setting

\begin{equation}
x = - e^{-\eta}, \quad y =  e^{\eta}, \quad z =- ( e^{\eta} - e^{-\eta})
\end{equation}

one obtains the skein relation

\begin{equation}
   e^{-\eta} \textsf{P} (\bar{L}_+) - e^{\eta}\textsf{P} (\bar{L}_-) = 
( e^{-\eta} - e^{\eta})
\textsf{P} (\bar{L}_0)
\end{equation}

   One can now exploit this relation along with ($5.21$) in well-known 
fashions to construct
invariant polynomials.( See also Ref.$14$ where a large number of 
sources are cited. ) We will
not present a full study of this aspect. Our aim has been to indicate 
the roles played by our
coefficients $e^{\pm\eta}$ ( an ingredient of our $\hat R$ ) in this 
context. This has been
achieved in our brief treatment.

\section{ From $\hat R$ to noncommutative spaces:}

\subsection{ Coordinates, differentials and mobile frames:}

  We implement well-known prescriptions [$15,16,17$] in the context of 
our class of $\hat R$. For
the $N^2\times N^2$ matrix $\hat R$ satisfying

\begin{equation}
(\hat R - I)(\hat R +e^{-2\eta}I) =0
\end{equation}

where for $\hat o(N),\hat p(N)$ respectively $$e^{\eta} +e^{-\eta} = 
[N\mp 1] \pm 1$$

let the coordinates ($x_1,x_2,...,x_N$) and the differentials 
($\xi_1,\xi_2,...,\xi_N$) be
ordered in in $N$-columns $x$ and $\xi$ respectively.

The prescriptions for the associated covariant differential 
geometries sitisfying the Leibnitz
rule [$15,16,17$] are

\begin{equation}
(1): \quad (\hat R - I)x\otimes x=0, \qquad x\otimes \xi =e^{2\eta} 
\hat R \xi\otimes x,\qquad
(\hat R +e^{-2\eta}I)\xi \otimes \xi =0
\end{equation}

\begin{equation}
(2): \quad (\hat R +e^{-2\eta} I)x\otimes x=0, \qquad x\otimes \xi =- 
\hat R \xi\otimes
x,\qquad (\hat R -I)\xi \otimes \xi =0
\end{equation}

  We concentrate below on ($6.2$). The  set ($6.3$) can be treated 
analogously , essentially
interchanging the roles of $x$ and $\xi$ ( except for the $(x,\xi)$ 
commutators ).

 From our previous definitions

$$ \hat R = I - e^{\mp \eta} P'_0$$

Hence from ($6.2$)

\begin{equation}
P'_0 x\otimes x =0
\end{equation}

The set of constraints ($6.4$) reduces to a ${\it single}$ one  due 
to the proportionality of the
nonzero rows of $P'_0$. This one  is easy to write down for all $N$ 
from ($1.17,1.18$). One
obtains for

\begin{equation}
\hat o(3): \qquad q^{-\frac {1}{2}}x_1x_3 +x_2x_2 + q^{\frac {1}{2}}x_3x_1 =0
\end{equation}

\begin{equation}
\hat o (4)): \qquad q^{-1}x_1x_4 +x_2x_3 +x_3x_2 + q x_4x_1 =0
\end{equation}

\begin{equation}
\hat p (4)): \qquad q^{-2}x_1x_4 +q^{-1}x_2x_3 -qx_3x_2 - q^2 x_4x_1 =0
\end{equation}

and so on.
  Consider now the constraints involving $\xi$ with the $\hat o (3)$ 
case as example.

Define

\begin{equation}
\Pi = (q^{-\frac {1}{2}}\xi_1x_3 +\xi_2x_2 + q^{\frac {1}{2}}\xi_3x_1)
\end{equation}

\begin{equation}
\Pi' = (q^{-\frac {1}{2}}\xi_1\xi_3 +\xi_2\xi_2 + q^{\frac {1}{2}}\xi_3\xi_1)
\end{equation}

  Now from ($6.2$) for $N=3$ with

$$e^{\eta}+e^{-\eta} = (q+1+q^{-1})$$

one has

\begin{equation}
x_i \xi_j = e^{2\eta}\xi_i x_j, \qquad (i+j \neq 4)
\end{equation}

\begin{equation}
x_1 \xi_3 = e^{2\eta}\xi_1 x_3 -e^{\eta}q^{-\frac {1}{2}}\Pi,\qquad 
x_2 \xi_2 = e^{2\eta}\xi_2
x_2 -e^{\eta}\Pi,\qquad x_3\xi_1=e^{2\eta}\xi_3 x_1 -e^{\eta}q^{\frac 
{1}{2}}\Pi
\end{equation}

\begin{equation}
\xi_i \xi_j =0, \qquad (i+j \neq 4)
\end{equation}

\begin{equation}
q^{-\frac {1}{2}}\xi_1 \xi_3 = y^{-1}q^{-1}\Pi',\quad \xi_2 \xi_2 =
y^{-1}\Pi', \quad q^{\frac {1}{2}}\xi_3 \xi_1 = y^{-1}q\Pi', \qquad y 
= (q+1+q^{-1})
\end{equation}

  Note the consistency of the sum of the equations ($6.13$). For $N >3$

\begin{equation}
\xi_i \xi_j =0, \qquad (i+j \neq (N+1))
\end{equation}

  Te coefficients of $(\Pi,\Pi')$ are now obtained from ($1.17,1.18$) 
they being proportional to
the nonzero elements in a row of $P'_0$. Also $\eta$ will now be as 
below ($6.1$).

  We now briefly consider the construction of mobile frames [$17,18$] 
( or "stehbeins" in the
terminology of the authors cited ). Let

\begin{equation}
\theta = \sum _{i} {\theta}_i {\xi}_i
\end{equation}

such that

\begin{equation}
[\theta,x_i] =0
\end{equation}

 From ($6.2,6.16$), remembering that in our conventions

$$\hat R = {\hat R}_{ij,kl} (ij)\otimes (kl)$$

and indicating all summations explicitly ( and using the $L^{\pm}$ of 
Secs.($3,4$) )
\begin{eqnarray}
\nonumber
(\sum _k \theta_k \xi_k)x_i&=& e^{-2\eta}\sum_{k,a,b}{\theta}_k 
{{\hat R}^{-1}}_{ka,ib}x_a\xi_b\\
\nonumber
&=& e^{-2\eta}\sum_{k,a,b} {\theta}_k ({{\hat R}^{-1}P})_{kb,ia} x_a \xi_b\\
\nonumber
&=& e^{-2\eta}\sum_{k,a,b} ({\theta}_k {L^-}_{kb,ia} x_a )\xi_b\\
&=& x_i(\sum_{b}{\theta}_b {\xi}_b)
\end{eqnarray}

Hence

\begin{equation}
x_i\theta_j = e^{-2\eta} \sum_{k,l} \theta_k L^{-}_{kj,il}x_l
\end{equation}

  Thus one obtains the commutators of ($x_i,\theta_j$). For $N=3$, setting

\begin{equation}
\theta_1x_1 +\theta_2x_2 +\theta_3x_1 = \tau
\end{equation}

$$ x_i \theta_i =  e^{-2\eta}\tau - e^{-\eta}\theta_{i'} x_{i'}, 
\qquad (i=1,2,3 )$$

$$ x_1\theta_2 = - e^{-\eta}q^{\frac {1}{2}}\theta_3 x_2,\qquad x_1\theta_3
=-e^{-\eta}\theta_3x_1$$
$$ x_2\theta_1 = - e^{-\eta} q^{\frac {1}{2}}\theta_2 x_3,\qquad x_2\theta_3
=-e^{-\eta} q^{-\frac {1}{2}} \theta_2x_1$$
\begin{equation}
  x_3\theta_1 = - e^{-\eta}\theta_1 x_3,\qquad x_3\theta_2
=-e^{-\eta} q^{-\frac {1}{2}}\theta_1x_2
\end{equation}

Generalizations for $N>3$ are obtained analogously.

  In Ref.$17$ different solutions of $\theta$ ( $3$ solutions for 
$So_q(3)$ ) are presented. They
involve a dilatation operator and inverses of coordinates, a radius 
$r$ being defined. We
consider no extensions of our algebras or such solutions for $\theta$ 
in the  present work.

\subsection{Towers of noncommutative ($x_1,..,x_N$) on a commutative 
($N-1$)- base space:}

  In eqn. ($1.3$) of Ref.$4$ it is pointed out that the mapping, 
implementing the transfer matrix
$t$,

\begin{equation}
\delta (x_i) = \sum_k t_{ik}\otimes x_k \equiv (tx)_i
\end{equation}

provides an iterative sequence of solutions. Since

\begin{equation}
\hat R (t \otimes I)(I\otimes t) = (t \otimes I)(I\otimes t)\hat R
\end{equation}

for  any polynomial of  $f(\hat R)$, if

        $$ f(\hat R)(x\otimes x) =0$$

\begin{equation}
  (t \otimes I)(I\otimes t) f(\hat R)(x\otimes x) =  f(\hat R)(t 
\otimes I)(I\otimes t)(x\otimes
x), \qquad = f(\hat R) (tx)\otimes (tx) =0
\end{equation}

  This is the mapping ($6.21$). It can evidently be iterated as

\begin{equation}
\delta ((t^n x)\otimes (t^n x)) =  (t^{n+1} x)\otimes (t^{n+1} x)
\end{equation}

We denote this as

$$ x^{(n)}  \rightarrow  x^{(n+1)} $$

  We note here a remarkable possibility, starting  again with $\hat 
o(3)$ as an example.

  Choose as the starting  point a ${\it commutative}$ solution of 
($6.5$) as follows:

  With parameters $ (a,b,c) \geq 0$ the following surface satisfies ($6.5$),

$$ x_1^{(0)} =a,\quad  x_3^{(0)} = -b,\quad x_2^{(0)}  = \pm \bigl( 
(q^{\frac {1}{2}}+q^{-\frac
{1}{2}}) ab \bigr)^{\frac {1}{2}}$$

\begin{equation}
  x_1^{(0)} =-a,\quad  x_3^{(0)} = b,\quad x_2^{(0)}  = \pm \bigl( 
(q^{\frac {1}{2}}+q^{-\frac
{1}{2}}) ab \bigr)^{\frac {1}{2}}
\end{equation}

  As ($a,b,c$) varies through real, non-negative values one  obtains a 
double cone whose
projections on the ($1,3$) plane covers the second and the  fourth 
quadrants. The vertices meet
at the origin. The projections of the contours $ x_2 =$ const. on the 
($1,3$) plane are
parabolas. The origin is invariant under $\delta$.

Now we implement $\delta$ as in ($6.21$) to obtain

\begin{equation}
x_i^{(1)} = (tx^{(0)})_i
\end{equation}

Consider for simplicity the $3\times 3$ fundamental $t$-matrices. 
These are given by ( compare
($4.1$) and see the reference cited above it)

\begin{equation}
t=t^+=  P \hat R=  R= P(\hat R P)P =L^+ _{21}
\end{equation}
and
\begin{equation}
t=t^-=  P \hat R^{-1}=  P(R^{-1}P)= ({\hat R}^{-1})_{21}  =L^- _{21}
\end{equation}

  We treat below $t^{\pm}$ together by setting correspondingly ( for $N=3)$ )

\begin{equation}
\lambda = \lambda_{\pm}= \lambda^{\pm 1}, \qquad (\lambda + \lambda ^{-1}+y =0)
\end{equation}

  For the fundamental rep. of $t$ our map gives

\begin{equation}
  x_1^{(1)} = \begin{vmatrix}  x_1^{(0)} &0 &0 \\ x_2^{(0)} &0 &0 \\ 
(1+q\lambda)x_3^{(0)}
&q^{\frac {1}{2}}\lambda x_2^{(0)}  &\lambda  x_1^{(0)}  \end{vmatrix}
\end{equation}

\begin{equation}
  x_2^{(1)}
  = \begin{vmatrix}  0 & x_1^{(0)} &0 \\ q^{\frac {1}{2}}\lambda 
x_3^{(0)} & (1+\lambda)x_2^{(0)}
&q^{-\frac {1}{2}}\lambda x_1^{(0)}  \\ 0 & x_3^{(0)}
&0 \end{vmatrix}
\end{equation}

\begin{equation}
  x_3^{(1)} = \begin{vmatrix}  \lambda x_3^{(0)} &q^{-\frac
{1}{2}}\lambda x_2^{(0)} &(1+q^{-1}\lambda)x_1^{(0)} \\ 0 &0 &x_2^{(0)}
\\ 0 &0  &  x_3^{(0)}  \end{vmatrix}
\end{equation}

The symmetries signalled above ($4.4$) reappear. Iteration now may proceed  as
$$  (x_i^{(0)}, x_i^{(1)}) \rightarrow  (x_i^{(n)}, x_i^{(n+1)}) $$

At each stage, given {\it only} ($6.5$) for $x_i^{(n)}$ ( and 
($6.29$) for $\lambda$ ), $x_i^{(n+1)}$
also satisfies ($6.5$) and hence
\begin{equation}
(\hat R-I) x^{(n+1)}\otimes  x^{(n+1)} = 0
\end{equation}

This has been verified explicitly. Moreover at each stage one can 
implement any chosen
representation of $t$ ( say, the $9\times 9$ rather than the $3 
\times 3$ ). Thus  one may
obtain varied sequences in the iterations. The illustration above is 
sufficient for our purpose.
  Let us compare this construction with a parallal possibility for 
$SO_q(3)$. We refer to the results
of Ex.$4.1.22$ of Ref.$17$. But for easier comparison with our 
results above we change the basis
from the circular components to our type as
  $$ (x_-,y,x_+) \rightarrow (x_1,x_2,x_3)$$

Now ($4.1.61$) of Ref.$17$ is

\begin{equation}
  x_1 x_2 = q x_2 x_1, \quad x_3 x_2 = q^{-1}x_2 x_3, \quad (x_3 x_1 - 
x_1 x_3) =
(q^{\frac{1}{2}} - q^{-\frac {1}{2}}) x_2^2
\end{equation}

There are now three constraints as compared to a single one for 
($6.5$). But one can choose the
{\it commutative} ($1,3$) plane (as compared to the  double cone 
before)  as the starting point by
setting

\begin{equation}
x^{(0)}_2 = 0,\quad x^{(0)}_1 x^{(0)}_3 =x^{(0)}_3 x^{(0)}_1
\end{equation}

This satisfies all the three constraints ($6.34$). Using the $3\times 
3$ $t$-matrix blocks for
$SO_q(3)$  (  and setting $\kappa = (q - q^{-1})$ ) one obtains

\begin{equation}
  x_1^{(n+1)} = \begin{vmatrix} q x_1^{(n)} &0 &0 \\0 &x_1^{(n)} &0 \\ 0
&0 &q^{-1} x_1^{(n)}  \end{vmatrix}
\end{equation}

\begin{equation}
  x_2^{(n+1)} = \begin{vmatrix}  x_2^{(n)} &\kappa x_1^{n} &0 \\0 
&x_2^{(n)} &q^{-\frac {1}{2}}
\kappa x_1^{n} \\ 0 &0 & x_2^{(n)}  \end{vmatrix}
\end{equation}
\begin{equation}
  x_3^{(n+1)} = \begin{vmatrix} q^{-1} x_3^{(n)} &q^{-\frac 
{1}{2}}\kappa x_2^{n} &\kappa
(1-q^{-1})x_1^{n}\\0 &x_3^{(n)} &\kappa x_2^{n} \\ 0 &0 & qx_3^{(n)} 
\end{vmatrix}
\end{equation}

Starting the iteration from ($6.35$) at each step $ x_i^{n}$ 
satisfies ($6.34$). As before one can
use more general realizations of the $t$-matrix at any step. One may 
note that for this case

         $$ (q=1) \rightarrow x_i^{(n+1)} = I_3 \otimes x_i^{(n)}$$

This is consistent with $\hat R = I$ for $q=1$ in the standard cases. 
But , as emphasized before,
there is no such triviality for our class for any value of $q$ 
(including $1$).  Note also that for
our class the matrices $x_i^{(1)}$ and the iterated ones are 
non-invertible. This is a  general
feature for our class.

   In the examples above one {\it starts} with a classical surface and 
iterating as above makes it
more and more "fuzzy" in this specific sense. This should be compared 
with the "fuzzy sphere" of
Ref.$17$ ( Sec.$7.2$) where one  starts fuzzy and a smooth surface is 
approached as a limit. One
moves in opposite senses in the two formalisms.

    So far we have studied the coordinate space ($x_i$) only. It must 
be noted carefully that one
cannot obtain a consistent nontrivial set $\xi_i^{(0)}$ 
corresponding to  $x_i^{(0)}$. This is
evident from ($6.10$) to ($6.13$) where $q$ is not restricted. So the 
whole covariant prescription
can be introduced at a noncommutative stage only. This however does 
not alter the fact that one can
build sequences of noncommutative ($x_i$) starting from a smooth surface.

   For $N>3$, with ($6.6$) and ($6.7$) as simplest examples, one has 
evidently more flexibility in
choosing $\xi_i^{(0)}$. Without going into details we indicate below 
for $\hat o(4)$  the stuctures
( valid more generally ) induced by our type of iterations.

$$ \hat o(4) : \qquad \lambda +\lambda^{-1} + (q^2 + 2 + q^{-2}) = 0 $$

\begin{equation}
x_1^{(n+1)} = \begin{vmatrix}  x_1^{(n)} &0 &0 &0 \\x_2^{(n)} &0&0&0 
\\x_3^{(n)}& 0&0&0 \\
(1+\lambda) x_4^{(n)}&q^{-1}\lambda x_3^{(n)}&q^{-1}\lambda 
x_2^{(n)}&q^{-2}\lambda x_1^{(n)}
   \end{vmatrix}
\end{equation}

\begin{equation}
x_2^{(n+1)} = \begin{vmatrix} 0& x_1^{(n)}  &0 &0 \\0&x_2^{(n)} &0&0 
\\q \lambda
  x_4^{(n)}&(1+\lambda) x_3^{(n)} &\lambda x_2^{(n)}&q^{-1}\lambda x_1^{(n)}
\\ 0& x_4^{(n)}&0&0 \end{vmatrix}
\end{equation}

\begin{equation}
x_3^{(n+1)} = \begin{vmatrix} 0&0& x_1^{(n)} &0 \\q \lambda
  x_4^{(n)}&\lambda x_3^{(n)} &(1+\lambda) x_2^{(n)}&q^{-1}\lambda 
x_1^{(n)} \\0&0
& x_3^{(n)}&0 \\ 0&0& x_4^{(n)}&0 \end{vmatrix}
\end{equation}

\begin{equation}
x_4^{(n+1)} = \begin{vmatrix} q^2\lambda x_4^{(n)} &q\lambda 
x_3^{(n)} &q\lambda x_2^{(n)}
&(1+\lambda) x_1^{(n)}\\ 0&0&0&x_2^{(n)}\\0&0&0&x_3^{(n)} \\ 
0&0&0&x_4^{(n)} \end{vmatrix}
\end{equation}

  The analogous structures for $\hat p(4)$ have also been obtained. 
The structures of the matrices
are quite similar to those for $\hat o(4)$. The $q$-dependence of the 
coefficients show typical
differences along with sign changes. They will not be reproduced 
here. In each case, for all $N$,
there is in each matrix one row and one column with nonzero elements. 
How they shift with with $i$
of $x_i^{(n)}$ should already be fairly apparent from the two 
preceding examples.

   \section {Remarks:}

   We have started to explore the contents of a class of braid 
matrices presented previously.
Unsurprisingly our study remains incomplete in all directions. We 
briefly indicate below
perspectives of further developments.

   For the  standard cases ( $q$-deformed unitary and orthogonal 
algebras ) we have studied
extensively elsewhere representations for $q$ a root of unity 
[19,20]. ( These two references cite
many other sources.) We introduced "the method of fractional parts" 
for this purpose. Here we have
noted (Sec.$2$) how certain roots of unity give triangularity ($\hat 
R^2 =I$) for different
dimensions. A study of our $L$-algebras for  such cases along similar 
lines can be of interest.

  Our other solutions for triangularity involve elliptic aand 
hyperelliptic $q$ with integer
(binomial) coefficients in  the defining equations. The roles of such 
functions  deserve further
explorations. Higher genus curves have appeared before [$21,22$] in 
the construction of
statistical models as solutions of star - triangle (or Yang-Baxter) 
relations. There such
curves  are necessary ingredients of the solutions. In our case the 
situation is quite differnt.
Our class of braid matrices have been obtained for {\it any} $q$. One 
can  even set $q=1$ and still
have interesting  solutions. Our special values of $q$ appear only 
when the additional constraint
of triangularity is imposed and depend on the dimension $N$.

  We have obtained  some important general features of  our 
$L$-algebras ( see  for  example the
eqns. from ($3.11$) to ($3.18$) ). But the explicit  study of 
realizations is limited to $3\times 3$
and $9\times 9$ ones for $ \hat o(3)$. This has already  shown the 
crucial role of the central,
group-like  elements we have constructed, thus achieving a principal 
goal. But a more general study
of the  $L$-algebras is desirable. Quadratic [$23$] and higher degree 
[$24$] homogenous algebras
have been studied from a functional point of  view yielding, for 
example, the Poincare series.
( More sources are cited in Ref.$24$.) The Poincare series for our 
algebra would show whether, and
if so what, irreducible representations interpolate the $N^{2^p} 
\times N^{2^p}$ dimensional
coproduct representations ( corresponding, as pointed out in Sec.$4$, 
to the central $\hat S_3$
elements proportional to $I_{N^{p}}$ ) obtained by iterating the 
coproduct prescription. We are
unable to answer this  question definitively at present, though 
attempts to realize intermediate
dimensional ones ( between $3\times 3$ and $9\times 9$ ) exploiting 
the  symmetries pointed out  in
Sec.$4$ seem to encounter obstructions. Our detailed study of  the 
$9\times 9$ case gives an idea of
the features to be  expected more generally.

  The Turaev construction for link invariants turns out to be 
elegantly adaptable to our case.
Systematic construction of invariant polynomials and possibility of 
generalizations to invariants of
$3$- manifolds will be  studied  elsewhere.

   Noncommutative geometries associated with our $\hat R$ have been 
presented indicating possible
constructions of "noncommutative towers" on classical base spaces of 
dimensions $< N $. Here again
a deeper study of the differential geometries remains to be done. 
Possible roles of our special
values of $q$ corresponding  to triangularity should be interesting 
to explore in this  context.

   It follows from (1.10) or ($3.41$) along with ($1.17$) that for $ 
\hat o(N)$ , real positive
$q$ and

                             $$ -\eta < \theta < 0$$

the elements of $\hat R (\theta)$ are all non-negetive ( either zero 
or real positive ). Hence such
an $\hat R (\theta)$ along with the  corresponding  transfer matrix 
$t(\theta)$ ( obtainable from
$\hat R (\theta)$ ) can furnish the basis of a multistate statistical 
model. The elements of $\hat R
(\theta)$ provide the Boltzmann weights. This class of models will be 
studied in a  following paper.


\begin{thebibliography}{}


\bibitem{1} A.Chakrabarti,J.Math.Phys.{\bf 44},5320 (2003)

\bibitem{2} A.Chakrabarti and R.Chakrabarti,J.Math.Phys.{\bf 44},785 (2003)

\bibitem{3} A.Chakrabarti,J.Math.Phys.{\bf 43},1589 (2002)

\bibitem {4} L.D.Faddeev,N.Yu.Reshetikhin and 
L.A.Takhtadzhyan,Leningrad Math.J.{\bf 1},193 (1990)

\bibitem {5} A.Klymik and K.Schmudgen, Quantum Groups and Their 
Representations, Springer (1997)

\bibitem {6} Tata Lectures on Theta II, D.Mumford,Birkhauser (1984) ( 
App.by H.Umemura)

\bibitem {7} R.Bruce King, Beyond the quartic equation, Birkhauser 
(1996) (Sec.8.3.-1)

\bibitem {8} B.Abdesselam,A.Chakrabarti and R.Chakrabarti, 
Mod.Phys.Lett. {\bf A13},779 (1998)

\bibitem {9} A.Chakrabarti and R.Chakrabarti, J.Phys.A:Math.Gen. {\bf 
33},1(2000)

\bibitem {10} A.P.Isaev, Sov.J.Part.Nucl.{\bf 26},501 (1995)

\bibitem {11} A.Chakrabarti,A nested sequence of projectors and 
corresponding braid matrices $\hat
R (\theta)$ : (I) Odd dimensions, math.QA/0401207 (App.B)

\bibitem{12} V.G.Turaev, Invent.Math.{\bf 92},527 (1988)

\bibitem {13} V.Chari and A.Pressley, Quantum Groups,C.U.P. (1994) 
(Secs.$15.1,15.2$)

\bibitem{14} F.Y.Wu,Rev.Mod.phys.{\bf 64},1099 (1992)

\bibitem{15} J.Wess and B.Zumino,Nucl.Phys.B (Proc.Suppl.){\bf 18B}, 302 (1990)

\bibitem{16} L.Hlavaty, J.Phys.A: Math.Gen. {\bf 25},485 (1992)

\bibitem{17} J.Madore, An Introduction to Noncommutative Differential 
Geometry,C.U.P. (1999)

\bibitem{18} B.L.Cerchiai,G.Fiore and J.Madore, Geometrical tools for 
quantum Euclidean spaces,
math.QA/0002007

\bibitem{19} B.Abdesselam,D.Arnaudon and A.Chakrabarti, 
J.Phys.A.Math.Gen.{\bf 28},5495 (1995)

\bibitem{20} B.Abdesselam,D.Arnaudon and A.Chakrabarti, 
J.Phys.A.Math.Gen.{\bf 28},3701 (1995)

\bibitem {21} H.Au-Yang,B.M.McCoy,J.Perk,S.Tang and M.L.Yan, 
Phys.Lett.{\bf A 123},219 (1987)

\bibitem {22} R.J.Baxter,J.H.H.Perk and H.Au-Yang, phys.Lett.{\bf 
A128},138 (1988)

\bibitem {23} Yu.I.Manin, Quantum Groups and Non-Commutative 
Geometry, CRM Univ. de Montreal (1988)

\bibitem {24} M.Dubois-Violette and T.Popov, Lett.Math.Phys.{\bf 61},159 (2002)

\end{thebibliography}
\end{document}